\documentclass[12pt,reqno]{amsart}
\usepackage{amsmath, amsthm, ams
symb, amsopn, amsfonts, amstext, enumerate, color, 
mathrsfs, url}
\usepackage[showonlyrefs]{mathtools}
\usepackage[shortlabels]{enumitem}
\usepackage{graphicx,color}
\usepackage{verbatim}
\usepackage[backref=page]{hyperref}  
\usepackage[square,comma,numbers,sort&compress]{natbib}

\usepackage[margin=1.2in]{geometry}

\makeatletter
\newcommand{\opnorm}{\@ifstar\@opnorms\@opnorm}
\newcommand{\@opnorms}[1]{%
 \left|\mkern-1.5mu\left|\mkern-1.5mu\left|
 #1
 \right|\mkern-1.5mu\right|\mkern-1.5mu\right|
}
\newcommand{\@opnorm}[2][]{%
 \mathopen{#1|\mkern-1.5mu#1|\mkern-1.5mu#1|}
 #2
 \mathclose{#1|\mkern-1.5mu#1|\mkern-1.5mu#1|}
}

\theoremstyle{plain}

\begin{document}


\theoremstyle{plain}
\newtheorem{theorem}{Theorem} [section]
\newtheorem{corollary}[theorem]{Corollary}
\newtheorem{lemma}[theorem]{Lemma}
\newtheorem{proposition}[theorem]{Proposition}


\theoremstyle{definition}
\newtheorem{definition}[theorem]{Definition}
\theoremstyle{remark}
\newtheorem{remark}[theorem]{Remark}

\numberwithin{theorem}{section}
\newtheorem*{theorem*}{Theorem}
\numberwithin{equation}{section}

\title[Estimates of stable solutions of Hardy-Henon equations in the ball]{A priori estimates of stable solutions of the general Hardy-Henon equation in the ball}

\thanks{The authors have been supported by the FEDER-MINECO Grant PID2021-122122NBI00 and by J. Andalucia (FQM-116).}

\author[J. Silverio Mart\'{i}nez-Baena]{J. Silverio Mart\'{\i}nez-Baena}
\address{J. Silverio Mart\'{i}nez-Baena\textsuperscript{1} ---
 Departamento de An\'{a}lisis Matem\'atico, Universidad de Granada, 18071 Granada, Spain.
}
\email{jsilverio@ugr.es}

\author[Salvador Villegas]{Salvador Villegas}
\address{Salvador Villegas\textsuperscript{2} ---
 Departamento de An\'alisis Matem\'atico, Universidad de Granada, 18071 Granada, Spain.
 }
\email{svillega@ugr.es}

\keywords{Hardy-Henon equation, regularity estimates, semi-stable solutions.
\\
\indent 2020 {\it Mathematics Subject Classification:}
35B45, 35B35, 35J61}
\date{}

\begin{abstract}
This paper is devoted to the study of semi-stable radial solutions
$u\in H^1(B_1)$ of $-\Delta u=\vert x\vert^\alpha f(u) \mbox{ in } B_1\setminus\lbrace0\rbrace$, where $f\in C^1(\mathbb{R})$ is a general nonlinearity, $\alpha>-2$  and $B_1$
is the unit ball of $\mathbb{R}^N$, $N>1$. We establish the boudness of such solutions for dimensions $2\leq N<10+4\alpha$ and sharp pointwise
estimates in the case $N\geq10+4\alpha$. 
In addition, we provide, for this range of dimensions, a large family of
semi-stable radially decreasing unbounded $H^1(B_1)$ solutions.
\end{abstract}

\maketitle


\section{Introduction}\label{section.intro}

This paper deals with the semi-stability of radial solutions $u\in H^1(B_1)$ of
 \begin{equation} \label{non-autonomous-equation}
-\Delta u  =  |x|^\alpha{f(u)} \qquad \mbox{in } B_1\setminus\lbrace0\rbrace\tag{P},
\end{equation}
where $x=(x_1,\dots,x_N)$ and $B_1$ the unit ball of $\mathbb{R}^N$ with measure $\omega_N/N$ for $N>1$. In the whole paper we assume $\alpha>-2$, and $f\in C^1(\mathbb{R})$. We deal with the regularity of semi-stable radially symmetric energy solutions $u\in H^1(B_1)$ of \eqref{non-autonomous-equation}. Denoting by $r:=|x|$ standard arguments show that these radial solutions $u=u(|x|)\equiv u(r)$ are continuous in $r\in (0,1]$ and can be viewed as classical solutions. We denote as usual $u_r:=\frac{<x,\nabla u>}{|x|}$ and constants by capital letters $C,K,C',K'\dots$ that may vary throughout the discussion and computations, with subscripts that determine the dependence of the constants on other parameters.

 The equation \eqref{non-autonomous-equation} corresponds to the Euler-Lagrange equation of the energy functional
 \begin{equation}
     E(u):=\int_{B_1}\left(|\nabla u|^2-|x|^\alpha F(u)\right)dx,
 \end{equation}
where $F(t)=\int_0^tf(s)ds$.
Since $f\in C^1(\mathbb{R})$ and $u$ is continuous, then $f'(u)$ is continuous. It is meaningful to call a solution $u\in H^1(B_1)$ semi-stable if the second variation of the energy is nonnegative definite, i.e.,
\begin{equation}\label{henon-stability-inequality}
    d^2E(\varphi,\varphi)=\int_{B_1}\left(|\nabla\varphi|^2-|x|^\alpha f'(u)\varphi^2\right) dx\geq 0,
\end{equation}
for all $\varphi\in C^1(B_1)$ with compact support in $B_1\setminus\lbrace 0 \rbrace$. 

Equation \eqref{non-autonomous-equation} can be regarded as a generalized version of the so called Hardy-Henon equation which corresponds to power-type nonlinearity $f(u)=|u|^{p-1}u$. That equation is in particular called the Lane-Emden equation in the autonomous case $\alpha=0$. Another nomenclature has been used for $f(u)=e^u$, the Henon-Gelfand problem. 

The study of this problem is preceded by the study of the so called Henon-equation whose origins goes back to the celebrated work of S. Chandrasekhar \cite[Chapter IV, pp. 87-88]{chandrasekhar1957stellar} on the stellar structure in polytropic equilibrium in astrophysics. 

The pionnering contribution in the autonomous case $\alpha =0$ is due to \cite{gidas-spruck} where the authors prove that for power-type nonlinearity $f(u)=|u|^{p-1}u$ and posed in the whole space $\mathbb{R}^N$, there is no positive solution for $p\in (1,p_S(N))$, where $p_S(N):=2^*-1=\frac{N+2}{N-2}$ if $N\leq 3$ and $p_S=+\infty$ if $N\leq 2$ is the classical critical Sobolev exponent.

For $p=p_{S}(N)$, the same equation is known to have (up to translation and rescaling) a unique positive solution, which is radial and explicit (see \cite{caffarelli-gidas-spruck}).

Let now $p_{JL}(N)>p_{S}(N)$ denote the so called Joseph-Lundgren exponent:
$$
p_{JL}(N)=\left\{
\begin{aligned}
+\infty&\quad\text{if $N\leq 10$},\\
\frac{(N-2)^2-{4N}+8\sqrt{N-1} }{(N-2)(N-10)}&\quad\text{if $N\ge 11$}.
\end{aligned}
\right.
$$
This exponent can be characterized as follows: for $p\ge p_{S}(n)$, the explicit function $u_{s}(x)=C_{p,N}\vert x\vert^{-\frac2{p-1}}$ for an appropiate constant $C_{p,N}$ depending only on $p$ and $N$, is a singular solution of the Lane-Emden equation, which is stable if and only if $p\geq p_{JL}(N)$.
It was proved in \cite{farina-classification} that Lane-Emden equation has no nontrivial finite Morse index solution whenever $1<p<p_{JL}(N)$, $p\neq p_{S}(N)$. Through the application of some blow-up analysis technique, such Liouville-type theorems imply interior regularity for solutions of a large class of semilinear elliptic equations: they are known to be equivalent to universal estimates for solutions of
\begin{equation} \label{general}
-Lu = f(x,u,\nabla u)\quad\text{in $\Omega$,}
\end{equation}
where $L$ is a uniformly elliptic operator with smooth coefficients, the nonlinearity $f$ scales like $\vert u\vert^{p-1}u$ for large values of $u$, and $\Omega$ is an open set of $\mathbb{R}^N$. For precise statements, see the work \cite{P-Q-S} in the subcritical setting, as well as its adaptation to the supercritical case by  \cite{davila-dupaigne-farina}.

For the Henon-Hardy equation (i.e. $\alpha\not=0$, $f(u)=|u|^{p-1}u$), there are also an extensive literature. In \cite{dancer-du-guo}, the authors show the following:

\begin{theorem}(Dancer, Du, Guo \cite{dancer-du-guo}). Let $u$ be a stable solution of the Henon-Hardy equation with $\alpha>-2$, $\Omega=\mathbb{R}^N$, and one of the following hypothesis is true: $ 2\leq N\leq 10+4\alpha$ and $ p>1$; or $ N>10+4\alpha$ and $ 1<p<p_{JL}(N,\alpha)$; where we have defined
    \[
    p_{JL}(N,\alpha)=\left\{
    \begin{array}{lr}
    +\infty\quad & \text{if } 2\leq N \leq 10+4\alpha,\\
    \frac{(N-2)^2-2(\alpha+2)(\alpha+N)+2\sqrt{(\alpha+2)^3(\alpha+2N-2)}}{(N-2)(N-4\alpha-10)}&\quad\text{if $N> 10+4\alpha$}.
    \end{array}
    \right.
    \]
Then $u\equiv0$. On the contrary, for $p\geq p_{JL}(N,\alpha)$ the equation admits a family of radially symmetric positive stable solutions.
\end{theorem}

 These results were later generalized in \cite{WangYe2012} for sign-changing solutions where the authors also proved that assuming $0\in\Omega$, no weak solution exists for $\alpha\leq-2$. As they remark in their proof, the result remains valid for any positive, convex and nondecreasing nonlinearity $f$. This indicates that the case $\alpha\leq -2$ carries additional pathological behavior, highlighting the necessity of avoiding that values for $\alpha$. For $\alpha>-2$, $\Omega=\mathbb{R}^N$ they also prove the following two results: 
\begin{enumerate}
    \item[$i)$] If $2\leq N<10+4\alpha$, there is no weak stable solution of the Gelfand-Henon equation.
    \item[$ii)$] If $2< N<10+4\alpha^-$ where $\alpha^-=\min\lbrace\alpha,0\rbrace$, then any weak solution of the Gelfand-Henon equation has infinite Morse index.
\end{enumerate}
\begin{remark}
    The dimension range $2\leq N<10+4\alpha$ is optimal, since for $N\geq 10+4\alpha$ $$u(x)=-(2+\alpha)\log|x|+\log[(2+\alpha)(N-2)],$$ is a radial stable weak solution to the Gelfand-Henon equation.
\end{remark}

For analogous recent results in the fractional/non local setting we refer to \cite{dupaigne-non-local}, \cite{FazlyHuYang2021}, \cite{hasegawa-ikoma-kawakami}, \cite{hasegawa-ikoma-kawakami-part-II}, \cite{HyderYang2021} and even with a Hardy potential term $|x|^{-2s}u$ on the right-hand side \cite{kim2022finite} or concerning systems \cite{dai2024liouville}, \cite{duong2024liouville} and references therein.

In our case, we are mainly concerned with the regularity of semi-stable radial solutions of the general Henon-Hardy equation, although we think that our results can be further developed in several directions involving extensions of what is known in the autonomous case. For instance, it would be interesting to investigate the existence and regularity of a critical solution of the equation $-\Delta u= \lambda |x|^\alpha f(u)$ with the appropriate set of hypothesis over $\alpha$ and $f$ in the spirit of \cite{brezis-vazquez}. The content is developed through a series of lemmas and propositions, arriving at the end to the proof of the main theorem stated at the beginning of the next section.

\section{Statement of the main result and steps of the proof}
We directly proceed to state the main result of the paper: a classification of the behavior of semi-stable radial solutions of \eqref{non-autonomous-equation} depending on the dimension.
\begin{theorem}\label{principal}

Let $N\geq 2$, $f\in C^1(\mathbb{R})$, $\alpha>-2$ and $u\in H^1(B_1)$ be a
semi-stable radial solution of (\ref{non-autonomous-equation}). Then there
exists a constant $C_{\alpha,N}$ depending only on $\alpha$ and $N$ such that:

\begin{enumerate}

\item[i)] If $N<10+4\alpha$, then $\Vert u\Vert_{L^\infty(B_1)}\leq C_{\alpha,N}
\Vert  u\Vert_{H^1 (B_1\setminus\overline{B_{1/2}})}$.

\

\item[ii)] If $N=10+4\alpha$, then $\vert u(r)\vert \leq C_{\alpha,10+4\alpha} \Vert
 u\Vert_{H^1 (B_1\setminus\overline{B_{1/2}})}\left( \vert \log r\vert +1\right) , \ \ \forall
r\in (0,1]$.

\

\item[iii)] If $N>10+4\alpha$, then $\displaystyle{\vert u(r)\vert \leq
C_{\alpha,N} \Vert  u\Vert_{H^1 (B_1\setminus\overline{B_{1/2}})}\, r^{\gamma(N,\alpha )}\, , \ \
\forall r\in (0,1]}$.
\end{enumerate}
where we have defined the exponent $\,\gamma(N,\alpha):=2-N/2+\alpha/2+\sqrt{(\alpha+2)(\alpha+2N-2)}/2$. 
\end{theorem}

\begin{remark}
    Observe that $\gamma(N,\alpha)=0$ if and only if $N=10+4\alpha$ and $\gamma(N,\alpha)<0$ if and only if $N>10+4\alpha$ so that the pointwise estimates $i)$ and $ii)$ do not give rise to boundedness of $u$.
\end{remark}

We shall divide the proof in several lemmas and propositions.
\begin{lemma}\label{lemma:key-non-autonomous} Let $N\geq 2$, $f\in C^1(\mathbb{R})$, $\alpha>-2$ and $u\in H^1(B_1)$ be a
semi-stable radial solution of \eqref{non-autonomous-equation}. Let $v\in C^{0,1}(0,1]$ such that $v(1)=0$. Then 
\begin{equation}
    \int_{r_0}^1 t^{N-1} u_r^2(t) \left(v'(t)^2 +\alpha \frac{v'(t)v(t)}{t}+\left(1-N-\alpha\frac{N}{2}\right)\frac{v(t)^2}{t^2} \right)dt\geq 0,
\end{equation}

\noindent for every $r_0\in (0,1)$.

\end{lemma}

\

\noindent {\bf Proof.} Assume, as a first step, that $v\in C^\infty (0,1)$ has compact support and $v\equiv 0$ in $(0,r_0]$ (later we will prove it for any $v\in C^{0,1}(0,1]$ with $v(1)=0$). Differentiating (\ref{non-autonomous-equation}) with respect to $r$ we obtain

$$(-\Delta u)_r=-\Delta u_r+\frac{N-1}{r^2}u_r=\alpha r^{\alpha -1}f(u)+r^\alpha f'(u)u_r.$$

Multiplying by the radial function $u_r v^2$ and integrating by parts yields

$$\int_{B_1}\left( \nabla u_r \nabla (u_r v^2)+\frac{N-1}{r^2}u_r^2 v^2 \right) dx=\int_{B_1}\left( \alpha r^{\alpha -1}f(u)+r^\alpha f'(u)u_r \right) u_r v^2 dx.$$

Therefore

\begin{equation}\label{hola}
\int_{B_1}\left( v^2 \vert\nabla u_r\vert^2+2u_r v\nabla u_r\nabla v+\frac{N-1}{r^2}u_r^2 v^2\right) dx= \int_{B_1}\left( \alpha r^{\alpha -1}f(u)u_r v^2+r^\alpha f'(u)u_r^2 v^2\right) dx.
\end{equation}

On the other hand, since $u$ is stable, we can consider the radial function $u_r v$ (which is a $C^1 (B_1)$ function with compact
support in $B_1\setminus\{ 0\}$) obtaining

\begin{equation}\label{laotra} \int_{B_1}  \vert \nabla (u_r v)\vert^2 dx\geq \int_{B_1} r^\alpha f'(u)(u_r v)^2 dx .
\end{equation}

Subtracting \eqref{hola} from \eqref{laotra} we can assert that

\begin{equation}\label{dfdf}
\int_{B_1}\left( u_r^2 v_r^2-\frac{N-1}{r^2}u_r^2 v^2\right)\geq-\alpha\int_{B_1}r^{\alpha -1}f(u)u_r v^2 dx.
\end{equation}

Let us expand this last term:

$$\int_{B_1}r^{\alpha -1}f(u)u_r v^2 dx=\int_{B_1}(-\Delta u)\frac{u_r v^2}{r}dx=\omega_N\int_{r_0}^1 t^{N-1}\left (-u_{rr}-\frac{N-1}{t}u_r \right)\frac{u_r v^2}{t}dt=$$

$$-\frac{\omega_N}{2}\int_{r_0}^1 \left(t^{2N-2}u_r^2\right)' t^{-N} v^2 dt=\frac{\omega_N}{2}\int_{r_0}^1 t^{2N-2}u_r^2 (t^{-N} v^2)' dt=$$

$$\omega_N\int_{r_0}^1 t^{N-1}\left( \frac{-N}{2}\frac{u_r^2 v^2}{t^2}+\frac{u_r^2 v v'}{t}\right)dt.$$

Combining this and \eqref{dfdf} we obtain 

$$\omega_N \int_{r_0}^1 t^{N-1}\left( u_r^2 v'^2-\frac{N-1}{t^2}u_r^2 v^2\right) dt\geq -\alpha\ \omega_N \int_{r_0}^1 t^{N-1}\left( \frac{-N}{2}\frac{u_r^2 v^2}{t^2}+\frac{u_r^2 v v'}{t}\right)dt,$$

\noindent which is the desired conclusion.

In fact, by standard density arguments, the above formula is also true if we consider a $C^{0,1}(0,1]$ function vanishing in $(0,r_0] \cup \{ 1\}$.

Now, for the sake of clarity of the exposition let us denote
\begin{equation}
   I(a,b;v):=\int_{a}^b t^{N-1} u_r^2(t) \left(v'(t)^2 +\alpha \frac{v'(t)v(t)}{t}+\left(1-N-\alpha\frac{N}{2}\right)\frac{v(t)^2}{t^2} \right)dt.
\end{equation}
For an arbitrary $v\in C^{0,1}(0,1]$ vanishing at $t=1$, we define a radial truncated function
\begin{equation}
    \bar{v}_\varepsilon(t)=
\begin{cases}
0&\text{if }0< t<\varepsilon,\\
\frac{v(r_0)}{r_0-\varepsilon}(t-\varepsilon)&\text{if }\varepsilon\leq t\leq r_0,\\
v(t)&\text{if }r_0<t\leq 1,
\end{cases}
\end{equation}
for any $0<\varepsilon<r_0$. Now, applying the first step of the lemma to the function $\bar{v}_\varepsilon$ we have that $I(\varepsilon,1;\bar{v}_\varepsilon)\geq 0$ and thus
\begin{align*}
    I(r_0,1; v)&\geq -I(\varepsilon,r_0;\bar{v}_\varepsilon)\\
    &=-\left(\frac{v(r_0)}{r_0-\varepsilon}\right)^2\int_{\varepsilon}^{r_0}t^{N-1}u_r(t)^2\left[1+\alpha\frac{t-\varepsilon}{t}+\left(1-N-\alpha\frac{N}{2}\right)\frac{(t-\varepsilon)^2}{t^2}\right]dt.
\end{align*}
Observe that since $u\in H^1(B_1)$ the whole right hand side is uniformly bounded for $\epsilon\in(0,r_0)$. Therefore 
\begin{align*}
    I(r_0,1; v)\geq -\lim_{\varepsilon\to 0}I(\varepsilon,r_0;\bar{v}_\varepsilon)
    =-\left(\frac{v(r_0)}{r_0}\right)^2(2+\alpha)\left(1-\frac{N}{2}\right)\int_{0}^{r_0}t^{N-1}u_r(t)^2dt
    \geq 0,
\end{align*}
as we claimed. \qed

\begin{remark}
    In the last step, it is shown that the assumption $u\in H^1(B_1)$ is essential. Indeed, it is well known (see for instance \cite{brezis-vazquez}), that there exists semi-stable weak solutions to $-\Delta u =C_{N,q}(1+u)^{\frac{q-2}{q}}$ in $B_1$ of the form $u(x)=|x|^q-1$ such that $u\notin H^1(B_1)$ in the range $q \in \left( -\frac{N}{2} + 2 - \sqrt{N - 1},\; -\frac{N}{2} + 1 \right]$, $N\geq 3$.
\end{remark}

\begin{proposition}\label{prop:non-zero-u_r}
    Let $N\geq 2$, $f\in C^1(\mathbb{R})$, $\alpha>-2$ and $u\in H^1(B_1)$ be a non constant semi-stable radial solution of \eqref{non-autonomous-equation}. Then $u_r\not=0$ in $(0,1]$.
\end{proposition}

\noindent {\bf Proof.} Assume by contradiction that there exists $r_1\in (0,1]$ 
 such that $u_r(r_1)=0$. From $u\in H^1(B_1)$, $f\in C^1(\mathbb{R})$ and the radial symmetry we know that $u(r)\in C^3((0,1])$. Thus, in particular there is a constant $C_{r_1}>0$ such that $|u_r(r)|\leq C_{r_1}|r-r_1|$ for any $r\in[r_1/2,r_1]$.
On the other hand, we take the following test function

\begin{equation*}
    v(t)=
\begin{cases}
\frac{t}{r_1-\varepsilon}&\text{if }0< t<r_1-\varepsilon,\\
\frac{r_1-t}{\varepsilon}&\text{if } r_1-\varepsilon\leq t\leq r_1,\\
0&\text{if }r_1<t\leq 1,
\end{cases}
\end{equation*}
with $0<\varepsilon<r_1/2$. Applying Lemma \ref{lemma:key-non-autonomous} we have that $-I(r_0,r_1-\varepsilon;v)\leq I(r_1-\varepsilon,r_1;v)$ for any $0<r_0<r_1-\varepsilon$. Therefore
\begin{align*}
    -\frac{1}{(r_1-\varepsilon)^2}&\int_{r_0}^{r_1-\varepsilon}t^{N-1} u_r^2(t)(2+\alpha)\left(1-\frac{N}{2}\right) dt\\ 
    &\leq\int_{r_1-\epsilon}^{r_1}\frac{t^{N-1} u_r^2(t)}{\varepsilon^2}\left[1-\alpha  \frac{(r_1-t)}{t}+\left(1-N-\alpha\frac{N}{2}\right)\frac{\left(r_1-t\right)^2}{t^2}\right]dt\\
    &\leq K \int_{r_1-\epsilon}^{r_1}\frac{ u_r^2(t)}{\varepsilon^2}dt\leq K C_{r_1}^{\,2}\int_{r_1-\epsilon}^{r_1}\frac{(t-r_1)^2}{\varepsilon^2}dt=\frac{K C_{r_1}^2}{3}\,\varepsilon.
\end{align*}
Taking limit as $\varepsilon\to 0$ we have that
\begin{equation*}
    -\frac{1}{r_1^2}(2+\alpha)\left(1-\frac{N}{2}\right)\int_{r_0}^{r_1}t^{N-1} u_r^2(t)  dt\leq0,
\end{equation*}
that implies that $u_r\equiv0$ in $[r_0,r_1]$ if $N>2$. By the uniqueness of the corresponding Cauchy problem, $u_r\equiv0$ in $(0,1]$ which gives the contradiction for $N>2$. Finally, if $N=2$ we change the test function by
\begin{equation*}
    v(t)=
\begin{cases}
\left(\frac{t}{r_1-\varepsilon}\right)^\beta&\text{if }0< t<r_1-\varepsilon,\\
\frac{r_1-t}{\varepsilon}&\text{if } r_1-\varepsilon\leq t\leq r_1,\\
0&\text{if }r_1<t\leq 1,
\end{cases}
\end{equation*}
with $\beta\in\mathbb{R}$ to be chosen later. With similar arguments we can conclude
\begin{equation*}
    -\frac{1}{r_1^{2\beta}}(\beta-1)\left(\beta+\alpha+1\right)\int_{r_0}^{r_1}t^{N-1+2\beta-2} u_r^2(t)  dt\leq0.
\end{equation*}
Chosing any $\beta\in(-1-\alpha,1)$, we obtain $u_r\equiv0$ in $[r_0,r_1]$ which is again a contradiction.\qed
\begin{lemma}\label{lemma:bound-of-u_r}
    Let $N\geq 2$, $f\in C^1(\mathbb{R})$, $\alpha>-2$ and $u\in H^1(B_1)$ be a non constant semi-stable radial solution of \eqref{non-autonomous-equation}. Then
    \begin{equation*}
        \int_{r/2}^ru_r^2\,dt\leq K_{\alpha,N}||\nabla u||^2_{L^2(B_1\setminus \overline{B_{1/2}})}\,r^{3-N+\alpha+\sqrt{(2+\alpha)(2N-2+\alpha)}},
    \end{equation*}
    for any $r\in (0,1]$.
\end{lemma}
\noindent {\bf Proof.} We first consider the case $0<r\leq 1/2$ and define
\begin{equation*}
    v(t)=
\begin{cases}
r^{s-1}t&\text{if }0< t<r,\\
t^s&\text{if } r\leq t\leq 1/2,\\
2^{1-s}(1-t)&\text{if }1/2<t\leq 1
\end{cases}
\end{equation*}
with $s\in\mathbb{R}$ to be determined. By Lemma \ref{lemma:key-non-autonomous} applied for $r_0=r/2$, we have
\begin{align*}
    0\leq r^{2s-2}&\int_{r/2}^r t^{N-1}u_r^2(t)\,(2+\alpha)\left(1-\frac{N}{2}\right)dt\\
    &+\int_{r}^{1/2}t^{N-3+2s}u_r^2(t)\left(s^2+\alpha s+1-N-\alpha\frac{N}{2}\right)dt\\
    &+\int_{1/2}^1t^{N-1}u_r^2(t)\,2^{2-2s}\left[1-\alpha \frac{1-t}{t}+\left(1-N-\alpha\frac{N}{2}\right)\left(\frac{1-t}{t}\right)^2\right]dt.
\end{align*}
Choosing $s_\alpha=-\alpha/2-\sqrt{(2+\alpha)(2N-2+\alpha)}/2$, the second term is indeed zero. Therefore, we deduce
\begin{align*}
     -r^{2s_\alpha-2}&\int_{r/2}^r t^{N-1}u_r^2(t)(2+\alpha)\left(1-\frac{N}{2}\right)dt\\
     &\leq \int_{1/2}^1t^{N-1}u_r^2(t)\,2^{2-2s_\alpha}\left[1-\alpha \frac{1-t}{t}+\left(1-N-\alpha\frac{N}{2}\right)\left(\frac{1-t}{t}\right)^2\right]dt.
\end{align*}
As a result
\begin{align*}
    - r^{2s_\alpha-2}\int_{r/2}^r & \left(\frac{r}{2}\right)^{N-1}u_r^2(t)(2+\alpha)\left(1-\frac{N}{2}\right)dt\\
    &\leq r^{2s_\alpha-2}\int_{r/2}^r t^{N-1}u_r^2(t)(2+\alpha)\left(1-\frac{N}{2}\right)dt\\
     &\leq \int_{1/2}^1t^{N-1}u_r^2(t)\,2^{2-2s_\alpha}\left[1-\alpha \frac{1-t}{t}+\left(1-N-\alpha\frac{N}{2}\right)\left(\frac{1-t}{t}\right)^2\right]dt\\
     &\leq C\,\int_{1/2}^1t^{N-1}u_r^2(t)dt,
\end{align*}
Finally, we obtain
\begin{align*}
    -\,  \frac{(2+\alpha)\left(1-\frac{N}{2}\right)}{2^{N-1}}r^{2s_\alpha-3+N}\int_{r/2}^r & u_r^2(t)dt\leq \frac{C}{\omega_N}||\nabla u||^2_{L^2(B_1\setminus\overline{B_{1/2}})}
\end{align*}
which is the desired conclusion if $N>2$.  If $N=2$, we change the test function by  
\begin{equation*}
    v(t)=
\begin{cases}
r^{-1-\alpha-\beta}t^\beta&\text{if }0< t<r,\\
t^{-1-\alpha}&\text{if } r\leq t\leq 1/2,\\
2^{2+\alpha}(1-t)&\text{if }1/2<t\leq 1.
\end{cases}
\end{equation*}
with $\beta\in (-1-\alpha,1)$ as in the previous proposition. We have
\begin{align*}
    -{r^{-2-2\alpha-2\beta}}\int_{r/2}^r & t^{-1+2\beta}\,u_r^2(t)(\beta-1)(\beta+\alpha+1)\,dt\\
    &\leq\int_{1/2}^1t\,u_r^2(t)\,2^{4+2\alpha}\left[1-\alpha \frac{1-t}{t}-\left(1+\alpha\right)\left(\frac{1-t}{t}\right)^2\right]dt,
    \end{align*}
that gives the desired conclusion once the same argument as above is employed for the last integral. 

On the other hand, let us call $R_{\alpha,N}$ a constant such that $0<R_{\alpha,N}<r^{3-N-2s_\alpha}$  for all $r\in[1/2,1]$. Hence, applying the above conclusion for $r=1/2$
\begin{align*}
    \int_{r/2}^r u_r^2(t)dt & \leq \int_{1/4}^{1/2} u_r^2(t)dt+\int_{1/2}^1u_r^2(t)dt\\
    &\leq K_{\alpha,N} \left(\frac{1}{2}\right)^{3-N-2s_\alpha}||\nabla u||^2_{L^2(B_1\setminus \overline{B_{1/2}})}+\frac{1}{\omega_N}||\nabla u||^2_{L^2(B_1\setminus \overline{B_{1/2}})}\\
    &= \left( K_{\alpha,N} \left(\frac{1}{2}\right)^{3-N-2s_\alpha}+\frac{1}{\omega_N}\right)\,\,||\nabla u||^2_{L^2(B_1\setminus \overline{B_{1/2}})}\\
    &\leq \frac{1}{R_{\alpha, N}}\left[K_{\alpha,N}\left(\frac{1}{2}\right)^{3-N-2s_\alpha}+\frac{1}{\omega_N}\right]||\nabla u||^2_{L^2(B_1\setminus \overline{B_{1/2}})}r^{3-N-2s_\alpha},
\end{align*}
giving rise to the same conclusion again.\qed

\begin{proposition}
    Let $N\geq 2$, $f\in C^1(\mathbb{R})$, $\alpha>-2$ and $u\in H^1(B_1)$ be a non constant semi-stable radial solution of \eqref{non-autonomous-equation}. Then
    \begin{align*}
        |u(r)-u(r/2)|\leq K'_{\alpha,N}||\nabla u||_{L^2(B_1\setminus \overline{B_{1/2}})}\,r^{2-N/2+\alpha/2+\sqrt{(2+\alpha)(2N-2+\alpha)}/2},
    \end{align*}
    for all $r\in(0,1]$.
\end{proposition}

\noindent {\bf Proof.} 
Fix $r\in(0,1/2)$, using Cauchy-Schwartz inequality and Lemma \ref{lemma:bound-of-u_r}, we have the following chain of inequalities:
\begin{align*}
    |u(r)-u(r/2)|= &\left|\int_{r/2}^ru_r(t)\,dt\right|\leq \int_{r/2}^r\left|u_r(t)\right|\,dt \leq\left(\int_{r/2}^ru_r^2(t)dt\right)^{1/2}\left(\int_{r/2}^rdt\right)^{1/2}\\
     &\leq K_{\alpha,N}||\nabla u||_{L^2(B_1\setminus \overline{B_{1/2}})}r^\frac{{3-N+\alpha+\sqrt{(2+\alpha)(2N-2+\alpha)}}}{2}\left(\frac{r}{2}\right)^{\frac{1}{2}}\\
     &\leq K'_{\alpha,N}||\nabla u||^2_{L^2(B_1\setminus \overline{B_{1/2}})}r^\frac{{4-N+\alpha+\sqrt{(2+\alpha)(2N-2+\alpha)}}}{2}.\ \ \ \ \ \ \ \ \ \ \ \ \ \ \ \ \ \ \ \ \ \ \ \ \ \   \qed
\end{align*}
\textbf{Proof of Theorem \ref{principal}.} 
Following the reasoning of \cite{villegas}, let $r\in(0,1]$, there exists $r_1\in(1/2,1]$ such that $r=r_1/2^{m-1}$ for some $m\in\mathbb{N}$. By Sobolev embedding in dimension one we observe that $u(r_1)\leq||u||_{L^\infty(B_1\setminus B_{1/2})}\leq \mathfrak{s}_N||u||_{H^1{(B_1\setminus \overline{B_{1/2}})}}$. Therefore
\begin{align}\label{eq:u(r)-bound}\nonumber
    |u(r)|\leq & |u(r_1)-u(r)|+|u(r_1)|\leq \sum_{k=1}^{m-1}\left|u\left(\frac{r_1}{2^{k-1}}\right)-u\left(\frac{r_1}{2^{k}}\right)\right|+|u(r_1)|\\\nonumber
    \leq & K'_{\alpha,N}||\nabla u||_{L^2(B_1\setminus \overline{B_{1/2}})}\sum_{k=1}^{m-1}\left(\frac{r_1}{2^{k-1}}\right)^{{2-\frac{N}{2}+\frac{\alpha+\sqrt{(2+\alpha)(2N-2+\alpha)}}{2}}}+\mathfrak{s}_{N}||u||_{H^1{(B_1\setminus \overline{B_{1/2}})}}\\
    \leq & \left(K''_{\alpha,N}\sum_{k=1}^{m-1}\left(\frac{r_1}{2^{k-1}}\right)^{{2-\frac{N}{2}+\frac{\alpha+\sqrt{(2+\alpha)(2N-2+\alpha)}}{2}}}+\mathfrak{s}_{N}\right)||u||_{H^1{(B_1\setminus \overline{B_{1/2}})}}.
\end{align}

$\bullet$ If $2\leq N<10+4\alpha$, we have $2-N/2+\alpha/2+\sqrt{(2+\alpha)(2N-2+\alpha)}/2>0$. Then
\begin{align*}
    \sum_{i=1}^{m-1}\left(\frac{r_1}{2^{i-1}}\right)^{2-\frac{N}{2}+\alpha/2+\sqrt{(2+\alpha)(2N-2+\alpha)}/2}\leq
\sum_{i=1}^{\infty}\left(\frac{1}{2^{i-1}}\right)^{2-\frac{N}{2}+\alpha/2+\sqrt{(2+\alpha)(2N-2+\alpha)}/2},
\end{align*}

\noindent which is a convergent series. Hence, equation \eqref{eq:u(r)-bound} implies statement i) of the theorem.

\

$\bullet$ If $N=10+4\alpha$, we have $2-N/2+\alpha/2+\sqrt{(2+\alpha)(2N-2+\alpha)}/2=0$. From \eqref{eq:u(r)-bound} we obtain
\begin{align*}
 \vert u(r)\vert &\leq \left( K'_{\alpha,N}
(m-1)+\mathfrak{s}_N \right) \Vert u\Vert_{H^1(B_1\setminus
\overline{B_{1/2}})}\\
&=\left( K'_{\alpha,N} \left( \frac{\log r_1-\log
r}{\log 2}\right) + \mathfrak{s}_N \right) \Vert
u\Vert_{H^1(B_1\setminus \overline{B_{1/2}})}\\ 
&\leq \left( \frac{K'_{\alpha,N}}{\log 2}+ \mathfrak{s}_N \right)
\left( \vert \log r\vert +1\right) \Vert u\Vert_{H^1(B_1\setminus
\overline{B_{1/2}})},
\end{align*}
\noindent which proves statement ii).

\

$\bullet$ Finally, if $N>10+4\alpha$, we have $2-N/2+\alpha/2+\sqrt{(2+\alpha)(2N-2+\alpha)}/2<0$. Then
\begin{align*}
    \sum_{i=1}^{m-1}\left(\frac{r_1}{2^{i-1}}\right)^{2-\frac{N}{2}+\frac{\alpha}{2}+\frac{\sqrt{(2+\alpha)(2N-2+\alpha)}}{2}}
    =\frac{r^{2-\frac{N}{2}+\frac{\alpha}{2}+\frac{\sqrt{(2+\alpha)(2N-2+\alpha)}}{2}}-r_1^{2-\frac{N}{2}+\frac{\alpha}{2}+\frac{\sqrt{(2+\alpha)(2N-2+\alpha)}}{2}}}{(1/2)^{2-\frac{N}{2}+\frac{\alpha}{2}+\frac{\sqrt{(2+\alpha)(2N-2+\alpha)}}{2}}-1}.
\end{align*}
From this and again \eqref{eq:u(r)-bound}, we deduce
\begin{align*}
    \vert u(r)\vert \leq \left(
\frac{K'_{\alpha,N}}{(1/2)^{2-\frac{N}{2}+\frac{\alpha}{2}+\frac{\sqrt{(2+\alpha)(2N-2+\alpha)}}{2}}-1}+\mathfrak{s}_N \right)
r^{2-\frac{N}{2}+\frac{\alpha}{2}+\frac{\sqrt{(2+\alpha)(2N-2+\alpha)}}{2}}\Vert u\Vert_{H^1(B_1\setminus
\overline{B_{1/2}})}\, ,
\end{align*}
which completes the proof. \qed
\section{Optimality of the Theorem \ref{principal}}
    These a priori estimates are indeed optimal as one can see with the following examples:
\begin{itemize}
    
     \item If $N=10+4\alpha$, then the function $u(x)=|\log |x||$ is a solution to the Henon-Gelfand equation
    \begin{equation*}
        -\Delta u = (N-2)|x|^{\alpha}e^{u(2+\alpha)}\qquad\text{in }B_1\setminus\lbrace0\rbrace,
    \end{equation*}
    such that $|x|^\alpha f'(u)=\frac{(N-2)^2}{4|x|^2}$. This implies that the stability inequality \eqref{henon-stability-inequality} attaining exactly the optimal Hardy constant and hence it is a stable solution.
    \\
    \item If $N>10+4\alpha$, recall that we have that for any $\alpha>-2$, $\gamma(\alpha,N)<0$. We consider the function $u(|x|)=|x|^{\gamma} -1$ for any 
    \begin{align}\label{eq:gamma-range-optimality-10+4alpha}
        \gamma(\alpha,N)\leq \gamma<0, 
    \end{align}
    which is the solution to the power-type equation
    \begin{equation*}
        -\Delta u =|x|^\alpha(-\gamma)(\gamma +N-2)(1+u)^{1+\frac{2+\alpha}{-\gamma}}\qquad\text{in }B_1\setminus\lbrace0\rbrace.
    \end{equation*}
    In this case 
    \begin{align}\label{stability-term-N-bigger-than-10+4alpha}
        |x|^\alpha f'(u)=\frac{(-\gamma+\alpha+2)(\gamma +N-2)}{|x|^2}.
    \end{align}
    Taking into account \eqref{eq:gamma-range-optimality-10+4alpha}, \eqref{stability-term-N-bigger-than-10+4alpha} and comparing again the stability inequality \eqref{henon-stability-inequality} with the hardy inequality, the constants satisfies:
    $$(-\gamma+\alpha+2)(\gamma +N-2)\leq \frac{(N-2)^2}{4},$$ 
    implying the stability of $u$.

\end{itemize}

\end{document}